\documentclass[letterpaper, 10 pt, conference]{ieeeconf}  

\IEEEoverridecommandlockouts                              

\usepackage{amsmath,amssymb}
\usepackage{amsfonts}
\usepackage{graphicx}
\usepackage{overpic}

\usepackage{color}
\usepackage{soul}
\usepackage[ruled,vlined]{algorithm2e}

\usepackage{theorem}

\newtheorem{corollary}{Corollary}
\newtheorem{lemma}{Lemma}

\newtheorem{proposition}{Proposition}

\newcommand{\Proof}{\noindent \textbf{Proof.}$\;\;$}
\newcommand{\BeginProof}{\noindent \textbf{Proof.}$\;\;$}
\newcommand{\EndProof}{\hfill{$\blacksquare$}\vspace{0.3cm}}

\newcommand{\vect}[1]{\ensuremath{\boldsymbol{\mathrm{#1}}}}
\newcommand {\matr}[2]{\left[\begin{array}{#1}#2\end{array}\right]}

\definecolor{seb}{rgb}{0.8,1,1}

\definecolor{kai}{rgb}{0.61, 0.87, 1.0}

\definecolor{todo}{rgb}{0.68, 1.0, 0.18} 

\title{\LARGE \bf
Recursive Feasibility of Stochastic Model Predictive Control with Mission-Wide Probabilistic Constraints
}
\author{Kai Wang and S\'ebastien Gros 
\thanks{All authors are with Department of Engineering Cybernetics, Norwegian University of Science and Technology (NTNU), 7491, Trondheim, Norway. {\tt\small\{kai.wang, sebastien.gros\}@ntnu.no} }
}

\begin{document}

\maketitle
\thispagestyle{empty}
\pagestyle{empty}

\begin{abstract}
This paper is concerned with solving chance-constrained finite-horizon optimal control problems, with a particular focus on the recursive feasibility issue of stochastic model predictive control (SMPC) in terms of mission-wide probability of safety (MWPS). MWPS assesses the probability that the entire state trajectory lies within the constraint set, and the objective of the SMPC controller is to ensure that it is no less than a threshold value. This differs from classic SMPC where the probability that the state lies in the constraint set is enforced independently at each time instant. Unlike robust MPC, where  strict recursive feasibility is satisfied by assuming that the uncertainty is supported by a compact set, the proposed concept of recursive feasibility for MWPS is based on the notion of remaining MWPSs, which is conserved in the expected value sense. We demonstrate the idea of mission-wide SMPC in the linear SMPC case by deploying a scenario-based algorithm. 
\end{abstract}

\section{Introduction}
Model predictive control (MPC) has been well established for dealing with complex constrained optimal control problems \cite{Rawlings2017model,Mayne2000constrained}. In the context of MPC, the system dynamics are required to be known and deterministic. In practice, the system uncertainty, including imprecise model parameters and process noises, is generally unavoidable. Because MPC does not take the uncertainty into account, constraints violations can occur. 

For applications where safety is critical, robust MPC (RMPC) strategies have been proposed to explicitly account for uncertainty. However, RMPC schemes can only handle bounded disturbances and the resulting control strategy can be conservative. To overcome these limitations, stochastic MPC (SMPC) methods have been developed to seek a trade-off between control performance and the risk of constraint violations using chance constraints.

There are mainly three forms of chance constraints proposed in the literature: individual chance constraint, stage-wise chance constraint \cite{Farina2016review} and mission-wide chance constraint \cite{Ono15}. E.g. in path planning for vehicles in the presence of obstacles, individual or stage-wise constraints restrict at every time instant the probability that the vehicle collides with an obstacle. In contrast, a mission-wide chance constraint directly restricts the probability of collision on the overall driving mission. A mission-wide chance constraint is arguably more meaningful than stage-wise constraints. Indeed, the former directly handles the risk of running a mission  \cite{Ono08,Ono15}, while the latter does it very indirectly.  However, stage-wise chance constraints are easier to handle than mission-wide constraints. Indeed, the latter handles probabilities over entire state trajectories, yielding very large probability spaces. More forms of chance constraints are discussed in, e.g., \cite[Section 2.2]{Farina2016review} and references therein.  

The current research on SMPC focuses on developping efficient methods for solving the underlying optimization problem, while recursive feasibility is less explored. Indeed, because of the (possibly unbounded) stochasticity, the recursive feasibility of SMPC typically holds in the probabilistic sense, making its analysis much more involved. Some results exist in specific contexts. When the system uncertainties have finite supports, recursive feasibility can be guaranteed using robust MPC \cite{Lorenzen2017Constraint}, at the cost of yielding very conservative control policies. For linear stochastic systems with infinite support, if the first two moments of the disturbance distribution are known,  constraint-tightening methods via the Chebyshev–Cantelli inequality are presented in \cite{Farina2015anApproach,Paulson2020joint,Yan2018Discounted}. Recursive feasibility is guaranteed using backup strategies when an infeasible optimization problem is encountered \cite{Farina2015anApproach,Paulson2020joint}, and using time-varying risk bound \cite{Yan2018Discounted}. The author in \cite{Ono2012recursive-a,Ono2012recursive-b} proposed SMPC algorithms that have a certain probability of remaining feasible if the initial condition is feasible. 
However, none of these methods tackle mission-wide probability of safety (MWPS), nor can provide a meaningful certificate of MWPS. In \cite{Pola2006Invariance}, the problem of maximizing the MWPS is expressed as a stochastic invariance problem and further developed into an optimal control problem, which is solvable via dynamic programming. Unfortunately, problems constraining the MWPS rather than maximizing it cannot necessarily be put in that simple form.

Guaranteeing recursive feasibility of a SMPC problem with MWPS constraints is an open problem, and this paper investigates a tentative solution. The main contributions of this paper is threefold. First, we show that if a policy is designed to achieve a certain MWPS, then the MWPS remaining until the end of the mission remains constant in the expected value sense. Second, we design a recursively feasible control scheme using shrinking horizon policies in the context of SMPC with MWPS guarantee. The proposed scheme treats directly the probability of running a mission successfully and therefore does not introduce artificial conservativeness. Third, an efficient scenario-based algorithm is proposed to deploy the idea in the linear case.

The paper is structured as follows. In Section \ref{sec:ProbState} we present the problem statement of SMPC with MWPS constraints, and its difference from the classical SMPC with stage-wise probabilistic constraint. Section \ref{sec:RecurFeas} details how the MWPS remains constant throughout the mission, and a recursively feasible policy design is discussed in Section \ref{sec:RecurFeas_MPC}. We demonstrate the idea in the linear SMPC case based on an efficient scenario-based algorithm in Section \ref{sec:LinearSMPC}. Finally, Section \ref{sec:conclusion} concludes the paper points to some future work.  

\textbf{Notation.$\,\,$} Boldface $\vect a$ (italic $a$) is a vector (scalar), and $\vect a_{0,\ldots,N}$ ($a_{0,\ldots,N}$) denotes a sequence of vectors (scalars). 
We use $\vect s_{0,\ldots,N}\in\mathbb{S}$ to denote that a state sequence $\vect s_{0,\ldots,N}$ lies in a constraint set $\mathbb{S}$ of the state space, i.e., $\vect s_{k} \in \mathbb S$ for all $k= 0,\ldots,N$. We denote $\mathbb I_{[a,b]}$ the set of integers in the interval $[a,b]\subseteq \mathbb R$.

\section{Problem Statement}
\label{sec:ProbState}
We consider a mission spanning a predefined horizon $N\in\mathbb N$ to be ``safe" if:
\begin{equation} \label{eq:safe}
	\vect s_{1,\ldots,N} \in \mathbb S  
\end{equation}
starting from some initial states $\vect s_0\in\mathbb S$. 
Here, $\vect s_k\in\mathbb R^{n}$ is the state at time step $k$, and $\mathbb S\subset\mathbb R^n$ is a set in the state space. We assume that the true stochastic system dynamics are given by:
\begin{equation}\label{eq:MC}
	\rho[\,\vect s_+\,|\, \vect s,\vect u\,]
\end{equation}
providing the probability density underlying transitions from a state-input pair $\vect s, \vect u$ to a new state $\vect s_+$. 
Throughout the paper, we assume that the states are known and continuous. Notice that the control community typically uses:
\begin{equation}\label{eq:dynmics}
\vect s_+ = \vect f(\vect s, \vect u, \vect w )
\end{equation}
to describe stochastic dynamics, in which $\vect w$ denotes the stochastic disturbances and $\vect f$ is generally a nonlinear function. The input $\vect u$ is given by a control policy sequence $$
	\vect\pi := \left\{\vect\pi_0,\ldots,\vect\pi_{N-1} \right\}$$
such that $$\vect u_k = \vect \pi_k(\vect s_k), \,\,  \forall k\in\mathbb I_{[0,N-1]}.$$

In general, guaranteeing the absolute safety as described in \eqref{eq:safe} yields very conservative control policies, or is even infeasible if the uncertainty is unbounded. Alternatively, for a given initial condition $\vect s_0\in\mathbb S$ and a policy sequence $\vect \pi$, we are interested in the Mission-Wide Probability of Safety (MWPS):
\begin{equation}
	\label{eq:MWPS}
	\mathbb{P}[ \vect s_{1,\ldots,N}\in\mathbb S\,|\, \vect s_{0},\vect \pi].
\end{equation}
The problem we are interested in is then to find a policy sequence solution of
\begin{subequations}\label{eq:GeneralOCP}
	\begin{align}
		\min_{\vect \pi}&\quad \mathbb E \left[M(\vect s_N)+\sum_{k=0}^{N-1}L(\vect s_k,\vect \pi_k(\vect s_k))\right] \label{eq:costs}\\
		\mathrm{s.t.}&\quad \mathbb{P}[\vect s_{1,\ldots,N} \in \mathbb S\,|\, \vect s_{0},\vect \pi]\geq S \label{eq:MWPC}\,,
	\end{align}
\end{subequations}  
where $S\in[0,1]$ is a predefined safety bound, the functions $L$ and $M$ are some given stage and terminal costs, and \eqref{eq:costs} is the expectation over the state trajectories resulting from $\vect s_0,\vect\pi$ and \eqref{eq:MC}.

In practice, calculating an optimal policy sequence for problem \eqref{eq:GeneralOCP} exactly is hardly possible, because it involves optimization over an infinite dimensional function space. To tackle this issue, we will be interested in using stochastic MPC formulations to generate policies that enforce the MWPS \eqref{eq:MWPS}, where a control input $\vect u_k$ is computed by solving an optimal control problem at every time step. In that context, a key concept will be the remaining MWPS at any time $k\in\mathbb I_{[1,N-1]}$ for a given state $\vect s_k$, defined as:
\begin{equation}\label{eq:cost:rest}
	\mathbb{P}[\vect s_{k+1,\ldots,N} \in \mathbb S\,|\, \vect s_{k}, \vect\pi]
\end{equation}
Here $\vect s_k$ is the outcome of a realization $\vect s_{1,\ldots,k}$ of the Markov Chain and its relationship to \eqref{eq:MWPS}.
 
A key observation is that a policy sequence $\vect\pi$ satisfying \eqref{eq:MWPC} does not yield any guarantee on \eqref{eq:cost:rest}. Indeed, an adversarial realization can e.g. bring the system into a state $\vect s_k$ for which the remaining MWPS is lower than $S$.

This observation entails that in the proposed context, the notion of recursive feasibility needs to be treated in a different way that is commonly done in robust MPC. We will detail this in Section \ref{sec:RecurFeas}. We briefly detail next the motivation for developing methods to treat MWPS constraints outside of the classical chance-constraint framework.

\subsection{Mission-wide constraints versus stage-wise constraints} \label{sec:MWPSvsSWPS}
In this section, we present our motivation for treating MWPS directly rather than via Stage-Wise Probability of Safety (SWPS). In particular, regardless of the desired MWPS level, enforcing it via SWPS becomes very conservative for long missions. SWPS problems seek policies that enforce constraints of the form:
\begin{equation}
	\label{eq:StageBound}
	\mathbb P\left[\vect s_k \in\mathbb S\,|\, \vect s_0,\,\vect \pi\right] \geq s_k\geq s ,\quad \forall\, k\in\mathbb I_{[1,N]}\,,
\end{equation} 
in which $1\geq s_k\geq s\geq 0$. One can then easily verify that the Boolean algebra and Booles's inequality entail that:
\begin{align}
	&\mathbb{P}[ \vect s_{1,\ldots,N}\in\mathbb S\,|\, \vect s_{0},\vect \pi] = 1 - \mathbb{P}\left[\left. \bigcup_{k=1}^{N} \vect s_{k}\notin\mathbb S\,\right|\, \vect s_{0},\vect \pi\right]  \nonumber\\
	&\geq 1 - \sum_{k=1}^{N} \mathbb{P}\left[ \vect s_{k}\notin\mathbb S\,|\, \vect s_{0},\vect \pi\right] \nonumber\\
	&= 1 -\sum_{k=1}^{N}\left(  1 -\mathbb{P}\left[ \vect s_{k}\in\mathbb S\,|\, \vect s_{0},\vect \pi\right] \right) \nonumber\\
	&\geq 1- N  + \sum_{k=1}^{N}s_k \geq 1 - N + Ns\,. \nonumber
\end{align}

\begin{figure}[htb]
	\centering
	\includegraphics[width=0.46\textwidth]{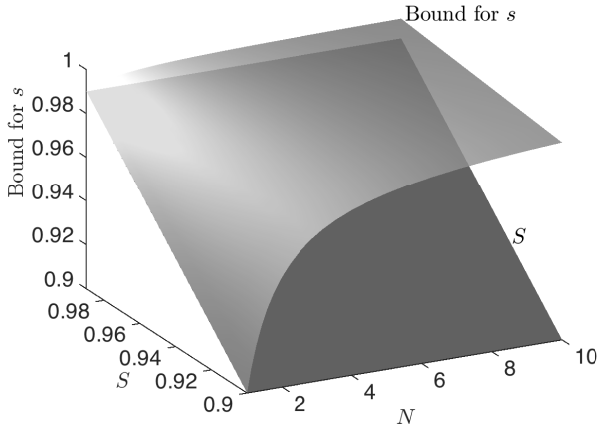}
	\caption{Illustration of bound \eqref{eq:MWPStoStageBound} for various $N$ and $S$. The curved manifold displays the bound \eqref{eq:MWPStoStageBound} for the SWPS such that a prescribed MWPS \eqref{eq:MWPS} holds.
	}
	\label{fig:parameters}
\end{figure}

To ensure the satisfaction of \eqref{eq:MWPC} via imposing \eqref{eq:StageBound}, requires the choice, $1- N  + Ns  \geq S$, i.e. a bound for $s$ can be derived as:
\begin{equation}\label{eq:MWPStoStageBound}
	s  \geq \frac{N - 1}{N} + \frac{S}{N}\,.
\end{equation}
Hence enforcing MWPS \eqref{eq:MWPS} via SWPS \eqref{eq:StageBound} requires selecting $s$ according to \eqref{eq:MWPStoStageBound}, which yields a bound $s$ close to one very fast as $N$ increases, see Fig. \ref{fig:parameters}, hence turning the SWPS constraints into hard constraints. While tighter bounds than \eqref{eq:MWPStoStageBound} can be derived\footnote{e.g. the Bonferroni inequalities allow one to refine the bound in \eqref{eq:MWPStoStageBound} by accounting for  some of the correlation between successive states}, treating MWPS via SWPS without introducing conservativeness is difficult. The intuitive reason behind this issue is that SWPS formulations neglect the time-correlation between the constraints violations, and as a result, it offers an incorrect representation of the risks incurred by a system over a mission. Bound \eqref{eq:MWPStoStageBound} corrects that, at the cost of introducing a high conservatism. By using a similar argument to what we developed above, risk-allocation technology proposed in \cite{OnoW08} optimizes the risk assigned  to each stage-wise constraint instead of using constant risk in \eqref{eq:StageBound}. This method leads a computationally expensive two-stage optimization problem and is still conservative as depicted in \cite[Fig. 1]{Farina2016review}.

\section{Relation between Remaining MWPS and Initial MWPS}\label{sec:RecurFeas}
Here, we show that the remaining MWPS is constant in the expected value sense. This offers a novel path for guaranteeing the recursive feasibility of MPC-like control schemes with MWPS constraints.

\begin{lemma} \label{lemma:FullPolicyFeasibility} 
	If the policy sequence $\vect\pi$ satisfies \eqref{eq:MWPS}, then
	\begin{align}
		\label{eq:Conservation:MWPS}
		\mathbb E_{\{\vect s_{1,\ldots,k}\in\mathbb S\,|\,\vect s_0,\vect\pi\}}\left[\mathbb{P}[ \vect s_{k+1,\ldots,N}\in\mathbb S\,|\, \vect s_{k},\vect \pi] \right]\geq S
	\end{align} 
	for all $k = 1,\ldots, N-1$, i.e. the remaining MWPS at time $k$ satisfies the constraint on the prescribed MWPS \eqref{eq:MWPS} in the expected value sense.
\end{lemma}

\Proof We observe that:
\begin{align*}
	&\mathbb{P}[ \vect s_{1,\ldots,N}\in\mathbb S\,|\, \vect s_{0},\vect \pi] \nonumber\\
	&=\int_{\mathbb S^{k}} \hspace{-0.1cm}\mathbb{P}[  \vect s_{1,\ldots,k}|\vect s_0,\vect\pi] \mathbb{P}[\vect s_{k+1,\ldots,N}\in\mathbb S\,| \vect s_{0,\ldots,k},\vect\pi] \mathrm d \vect s_{1}\ldots\mathrm d \vect s_{k}  \nonumber \\
	&=\int_{\mathbb S^{k}}\mathbb{P}[\vect s_{1,\ldots,k}\,|\,\vect s_0,\vect\pi]\mathbb{P}[ \vect s_{k+1,\ldots,N}\in\mathbb S\,|\, \vect s_{k},\vect\pi]  \,\mathrm d \vect s_{1}\ldots\mathrm d \vect s_{k} \nonumber\\
	&=\int_{\mathbb S} \mathbb{P}[ \vect s_{1,\ldots,k-1}\in\mathbb S\wedge\vect s_k|\vect s_0,\vect\pi] \mathbb{P} [\vect s_{k+1,\ldots,N}\in\mathbb S\,| \vect s_{k},\vect\pi] \mathrm d \vect s_k \nonumber\\
	&:=\,\mathbb E_{\{\vect s_{1,\ldots,k}\in\mathbb S\,|\,\vect s_0,\vect\pi\}}\left[\, \mathbb{P}[ \vect s_{k+1,\ldots,N}\in\mathbb S\,|\, \vect s_{k},\vect\pi]\, \right]\,.   
\end{align*}
Here, $\vect s_{1,\ldots,k}$ is a Markov Chain underlying \eqref{eq:MC}, and therefore a random variable in the high dimensional space $\left(\mathbb R^n\right)^k$. Given a policy sequence $\vect \pi$, the remaining MWPS at time $k$, i.e., the term inside $\mathbb E_{\{\vect s_{1,\ldots,k}\in\mathbb S\,|\,\vect s_0,\vect\pi\}}[\,\cdot\,]$, depends only on the random state $\vect s_k$, which is a particular dimension in the Markov Chain $\vect s_{1,\ldots,k}$. Hence, the last equation holds because here $\mathbb E_{\{\vect s_{1,\ldots,k}\in\mathbb S\,|\,\vect s_0,\vect\pi\}}[\,\cdot\,]$ is used to denote the expectation value of the remaining MWPS that is taken over all possible realizations of the random Markov Chain $\vect s_{1,\ldots,k}$ that remains in $\mathbb S$. \EndProof

Lemma \ref{lemma:FullPolicyFeasibility} entails that the MWPS is conserved in the expected value sense throughout the mission if a mission-wide policy sequence has been selected at the beginning of the mission. Result \eqref{eq:Conservation:MWPS} is arguably best interpreted in a frequentist framework. Indeed, even though a specific realization $\vect s_{0,\ldots,k}$ may be adversarial for the remaining MWPS, we observe that in average the MWPS remains unchanged throughout the mission. As a result, \eqref{eq:Conservation:MWPS} entails that if running missions under policy $\vect\pi$ designed according to \eqref{eq:GeneralOCP}, the resulting ratio of success will asymptotically be at least $S$. While this statement may appear tautological, it provides a basic concept of recursive feasibility that can be translated into constraints in a MPC framework to ensure that a prescribed MWPS is achieved. We detail this observation below.


\section{Recursive feasibility of MWPS with shrinking-horizon policies}\label{sec:RecurFeas_MPC}
In this section, we focus on solving the originally proposed mission-wide probability-constrained finite-horizon optimal control problem \eqref{eq:GeneralOCP} using shrinking-horizon policies that are updated as the mission progresses. The reason behind this is that the exact optimal policies for \eqref{eq:GeneralOCP} is difficult to compute in general. 

As a result, in practice, the policy sequence $\vect\pi$ is typically finitely parameterized,
	and hence restricted to a subset of the set of admissible policies.
	This introduces sub-optimality, and makes it useful to re-solve problem \eqref{eq:GeneralOCP} at every time instant $k$, according
	to the latest state realization $\vect s_k$. We then consider at every time $k$ the control policy sequence: 
\begin{equation}\label{eq:PredictedPolicy}
	\vect \pi^k = \left\{\vect \pi^k_{k},\ldots,\vect \pi^k_{N-1}\right\}
\end{equation}
lasting to the end of the mission. For the sake of brevity, we will work with a shrinking horizon extending to the end of the mission. The fixed, receding horizon shorter than the mission duration will be the object of our future work.

At every time instant $k\in\mathbb I_{[0,N-1]}$, for the corresponding state $\vect s_k$, we consider solving the following shrinking-horizon, mission-wide and chance-constrained problem:
\begin{subequations}\label{eq:OCPinSMPC}
	\begin{align}
		\min_{\vect \pi^k}&\quad \mathbb E \left[M(\vect s_{N})+\sum_{l=k}^{N-1}L(\vect s_{l},\vect \pi^k_{l}(\vect s_{l}))\right] \label{eq:costsinSMPC}\\
		\mathrm{s.t.}&\quad \mathbb{P}[\vect s_{k+1,\ldots,N} \in \mathbb S\,|\, \vect s_{k},\vect \pi^k]\geq S_k \label{eq:MWPS_1_in_SMPC}
	\end{align}
\end{subequations}  
to get a new policy sequence. Here, $S_k\in[0,1]$ is a varying risk-bound that will be specified later. Notice that while \eqref{eq:OCPinSMPC} computes an entire policy sequence $\vect \pi^k$ for the current state $\vect s_k$, only the first policy $\vect \pi^k_k$ of that sequence is used to generate the actual control action, as the policy sequence is recalculated at the next time instant $k+1$, in a classic MPC fashion. The inputs eventually applied to the closed-loop system will therefore read as:
\begin{equation}
	\vect u_k = \vect\pi^k_k\left(\vect s_k\right)\,, \quad \forall\,k\in\mathbb I_{[0,N-1]}
\end{equation}

In the context of mission-wide SMPC, we will consider the recursive feasibility issue of employing the policy sequence $\{\vect \pi^0_0,\ldots,\vect \pi^{N-1}_{N-1}\}$ resulting from extracting only the first policy $\vect \pi^k_k$ of the policy sequence $\vect \pi^k$ at every time step $k$, for all $k\in\mathbb I_{[0,N-1]}$. We show next that retaining recursive feasibility in the sense of \eqref{eq:Conservation:MWPS} requires only that the new policy sequence produces a remaining MWPS that is not worse than a discounted one achieved by the previous policy for the current state $\vect s_k$. We formalise this statement in the proposition below.

\begin{proposition} \label{Prop:InductivePolicies_1} 
	Assume that the initial policy sequence $\vect\pi^0 $ satisfies the MWPS cosntraint:
	\begin{equation}
		\label{eq:MWPS:Conservative_1}
		\mathbb{P}[ \vect s_{1,\ldots,N}\in\mathbb S\,|\, \vect s_{0},\vect \pi^0 ] \geq S_0 \geq S
	\end{equation}
	and that each policy sequence $\vect\pi^k$ is built under the constraint:
	\begin{equation}
		\label{eq:MWPS:nondecreasing:Conservative_1}
		\mathbb{P}[ \vect s_{k+1,\ldots,N}\in\mathbb S\,|\, \vect s_{k},\vect \pi^k ] \geq  S_k
	\end{equation}
	where
	\begin{equation*}
		S_{k} = \gamma_k \mathbb{P}[ \vect s_{k+1,\ldots,N}\in\mathbb S\,|\, \vect s_{k},\vect \pi^{k-1} ]
	\end{equation*}
	holds and 
	with $\gamma_k\in(0,1]$, for all $k\in\mathbb I_{[1,N-1]}$. Then the MWPS under $\vect u_k = \vect\pi_k^k\left(\vect s_k\right)$ and $k\in\mathbb I_{[0,N-1]}$ reads as:
	\begin{equation}
		\label{eq:InductivePolicies_1}
		\mathbb{P}\left[ \vect s_{1,\ldots,N}\in\mathbb S\,\middle |\, \vect s_{0},\{\vect \pi^0_0,\ldots,\vect \pi^{N-1}_{N-1}\}\right] \geq \prod_{k=1}^{N-1}\gamma_k S_0\,.
	\end{equation}
\end{proposition}
\Proof We will prove this by induction. Consider
\begin{align*}
	&\mathbb{P}\left[ \vect s_{1,\ldots,N}\in\mathbb S\,\middle |\, \vect s_{0},\{\vect \pi^0_0,\ldots,\vect \pi^{k}_{k},\ldots,\vect\pi^{k}_{N-1}\}\right] \nonumber\\
	&=\int_{\mathbb S}\mathbb{P}\left[ \vect s_{1,\ldots,k-1}\in\mathbb S\wedge \vect s_k\,|\, \vect s_{0},\{\vect \pi^0_0,\ldots,\vect \pi^{k-1}_{k-1}\}\right] \nonumber\\
	&\hspace{3.5cm}\cdot \mathbb{P}\left[ \vect s_{k+1,\ldots,N}\in\mathbb S\,\middle |\, \vect s_{k},\vect\pi^{k}\right] \,\mathrm d\vect s_k \nonumber \\
	&\geq\int_{\mathbb S}\mathbb{P}\left[ \vect s_{1,\ldots,k-1}\in\mathbb S\wedge \vect s_k\,|\, \vect s_{0},\{\vect \pi^0_0,\ldots,\vect \pi^{k-1}_{k-1}\}\right] \nonumber\\
	&\hspace{3cm} \cdot\gamma_k   \mathbb{P}\left[ \vect s_{k+1,\ldots,N}\in\mathbb S\,\middle |\, \vect s_{k},\vect \pi^{k-1} \right]\mathrm d\vect s_k \nonumber\\
	&=\gamma_k\mathbb{P}\left[ \vect s_{1,\ldots,N}\in\mathbb S\,\middle |\, \vect s_{0},\{\vect \pi^0_0,\ldots,\vect \pi^{k-1}_{k-1},\ldots,\vect\pi^{k-1}_{N-1}\}\right]\,, \nonumber  
\end{align*}
where the last equality holds because the last integral describes the MWPS associated to applying the policy sequence $\{\vect \pi^0_0,\ldots,\vect \pi^{k-1}_{k-1},\ldots,\vect\pi^{k-1}_{N-1}\}$. Hence an induction from
\begin{align}\nonumber
 \mathbb{P}[ \vect s_{1,\ldots,N}\in\mathbb S\,|\, \vect s_{0},\vect\pi^{0}] \geq S_0
\end{align}
yields \eqref{eq:InductivePolicies_1}.\EndProof

Let us introduce the following corollaries, showing the practical implications of Proposition \ref{Prop:InductivePolicies_1}:
\begin{corollary} \textbf{Guarantee of MWPS:}
	The choice:
	\begin{align} \label{eq:MPWS:design1_1}
		\prod_{k=1}^{N-1}\gamma_k S_0 = S
	\end{align}
	together with the policy update constraint \eqref{eq:MWPS:nondecreasing:Conservative_1} yields a sequence of policies $\{\vect \pi^0_0,\ldots,\vect \pi^{N-1}_{N-1}\}$ that satisfies the prescribed MWPS constraint \eqref{eq:MWPS:Conservative_1}.
\end{corollary}
\BeginProof The update constraint \eqref{eq:MWPS:nondecreasing:Conservative_1} ensures that \eqref{eq:InductivePolicies_1} holds. Condition \eqref{eq:MPWS:design1_1} imposed on the factors $\gamma_{ 1,\ldots,N-1}$ then ensures that \eqref{eq:MWPC} is satisfied.\EndProof

\begin{corollary}\textbf{Recursive Feasibility:} Constraint \eqref{eq:MWPS:nondecreasing:Conservative_1} is always feasible for any $\gamma \leq 1$
\end{corollary}
\BeginProof
We observe that \eqref{eq:MWPS:nondecreasing:Conservative_1} is feasible for $\vect\pi^k = \vect\pi^{k-1}$.
\EndProof

\section{A Scenario-based Linear SMPC approach with Mission-Wide Guarantees}
\label{sec:LinearSMPC}
In this section we deploy the mission-wide SMPC idea developed so far in the linear case. Let us consider that the stochastic dynamics \eqref{eq:dynmics} are explicitly given by:
\begin{align}\label{eq:LinearSystem}
	\vect s_{k+1} = A\vect s_k + B\vect u_k +  \vect w_k\,,
\end{align}
and that the safe set $\mathbb S$ is polytopic, i.e.
\begin{align}\label{eq:safeSet}
	\mathbb S = \left\{\,\vect s\,|\, C\vect s+\vect c \leq 0\right\}\,.
\end{align} 
Here we assume that the disturbances $\vect w_k$, $k\in\mathbb I_{[0,N-1]}$ are i.i.d., and zero-mean for the sake of notation convenience. 

At each time instants $k\in\mathbb{I}_{[0,N-1]}$, the predicted state $\vect s_t$ for all  $t =k,k+1,\ldots,N$, can be split into a nominal part and an stochastic error part, i.e., $\vect s_t = \bar{\vect s}_t + \vect e_t$. we consider the policy sequence $\vect\pi^k_t$ parameterized via $\bar{\vect u}_t,\, K$, given by:
$$\vect u_t = \vect\pi^k_t\left(\vect s_t\right) := \bar{\vect u}_t + K\vect e_t, \quad \forall t \in\mathbb{I}_{[k,N-1]}$$ where $K$ is a stabilizing feedback matrix for the nominal dynamics:
\begin{align}
	\bar{\vect s}_{t+1} = A\bar{\vect s}_t + B\bar{\vect u}_t\,\quad \bar{\vect s}_{k} = \vect s_k.
\end{align} 
The stochastic error dynamics are then given by:
\begin{align}
	\vect e_{t+1} = \left(A+BK\right)\vect e_{t}  +  \vect w_t,\quad \vect e_k = 0\,.
\end{align}

Our goal is to solve the following mission-wide probability constrained optimal control problem at every time instant $k$:
\begin{subequations}\label{eq:LinearSMPC:OCP}
	\begin{align}
		\min_{\bar{\vect u}_{k,\ldots,N-1}}&\,\, \mathbb E \left[\vect s_N^{\top} Q_N \vect s_N + \sum_{t=k}^{N-1}\left(\vect s_t^{\top} Q \vect s_t + \vect u_t^{\top} R \vect u_t\right)\right]\label{cost}\\
		\mathrm{s.t.\,\,\,\,\,}&\,\, \bar{\vect s}_{k}=\vect s_k\\
		&\,\, \bar{\vect s}_{t+1} = A\bar{\vect s}_t + B\bar{\vect u}_t,\hspace{1.3cm} \forall t\in\mathbb I_{[k,N-1]}\\
		&\,\, \vect e_{t+1} = \left(A+BK\right)\vect e_{t}  +  \vect w_t,\,\, \forall t\in\mathbb I_{[k,N-1]}\label{chance1}\\
		&\,\, \vect s_{t+1} = \bar{\vect s}_{t+1}+\vect e_{t+1},\hspace{1.2cm} \forall t\in\mathbb I_{[k,N-1]}\label{chance2}\\
		&\,\,\mathbb P [C\vect s_{t+1}+\vect c \leq 0,\,\, \forall t\in\mathbb I_{[k,N-1]}] \geq S_k\,.\label{chance}
	\end{align}
\end{subequations}   
Here, $Q, Q_N$ are semi-positive definite,  $R$ is positive definite,
and the value
\begin{equation} \label{Sk}
	S_k = \gamma_k\mathbb{P}[\vect s_{k+1,\ldots,N} \in \mathbb S\,|\,\vect s_{k}, \vect\pi^{k-1}]
\end{equation}
will be estimated at every time instant $k$ using Monte Carlo simulation based on the real, closed-loop state $\vect s_k$ and the previous policy sequence $\vect\pi^{k-1}$.

\subsection{An efficient scenario-based SMPC algorithm}
\vspace{0.2cm}
\noindent\textbf{Cost Function.} Since $\mathbb E [\vect s_t] = \bar{\vect s}_t$ and $\vect e_t$ is zero mean (since $\vect w_t$ is zero mean by assumption), the cost function \eqref{cost} can be written explicitly as 
\[
\bar{\vect s}_N^{\top} Q_N \bar{\vect s}_N + \sum_{t=k}^{N-1}\left(\bar{\vect s}_t^{\top} Q \bar{\vect s}_t + \bar{\vect u}_t^{\top} R \bar{\vect u}_t\right) + \sigma, 
\]
where $\sigma$ is a constant term that can be excluded from the cost function.

\vspace{0.2cm}
\noindent\textbf{Chance Constraint.} Substituting \eqref{chance1} and \eqref{chance2} into \eqref{chance}, the constraints can be rewritten as
\[
\mathbb P [C(\bar{\vect s}_{t+1} + \left(A+BK\right)\vect e_{t}  +  \vect w_t)+\vect c \leq 0, \forall t\in\mathbb I_{[k,N-1]}] \hspace{-0.2em}\geq\hspace{-0.2em} S_k 
\]
and further be written as
\begin{equation}\label{chance constraint}
\mathbb P [\underbrace{\mathcal C \mathcal A\vect w_{k,\ldots,N-1}^\top+[\vect c,\ldots,\vect c]^\top}_{:=H} \hspace{-0.2em}+ \,\mathcal C \bar{\vect s}_{k+1,\ldots,N}^\top\leq \hspace{-0.2em}\vect 0] \hspace{-0.2em}\geq\hspace{-0.2em} S_k 
\end{equation}
with matrix $\mathcal A$ obtained by condensing the dynamic \eqref{chance1} and matrix $\mathcal C$ being block-diagonal with $C$ as blocks.

\vspace{0.2cm}
\noindent\textbf{Scenario Approximation.} In general, providing a closed form for \eqref{chance constraint} is difficult. Fortunately, this problem can be handled efficiently with a scenario-based approach. Constraints \eqref{chance constraint} is replaced by a finite, sufficiently large number $N_k$ of deterministic constraints resulting from sampling the disturbance sequence $\vect w_{k,\ldots, N-1}$. For a given time instant $k$, we define the $i^{\mathrm{th}}$ sample for all $i\in\mathbb I_{[1,N_k]}$ as
$$
\vect w^{(i)}_{k,\ldots,N-1} := \{\vect w^{(i)}_k,\ldots,\vect w^{(i)}_{N-1}\}\,,
$$  
Hence, the chance constraint \eqref{chance constraint} can be converted to 
\begin{equation}\label{deterministic}
H^{(i)} + \,\mathcal C \bar{\vect s}_{k+1,\ldots,N}^\top\leq \vect 0, \quad \forall i\in\mathbb I_{[1,N_k]}\,,
\end{equation}
where
\[
H^{(i)}=\mathcal C \mathcal A(\vect w^{(i)}_{k,\ldots,N-1})^\top+[\vect c,\ldots,\vect c]^\top\,.
\]

In order to guarantee that \eqref{deterministic} approximates \eqref{chance constraint}  with a high probability $1-\beta$, where $\beta$ is typically set to be very small (e.g., $\beta=10^{-6}$), $N_k$ must satisfy the following inequality \cite{Calafiore2006}:
\begin{equation*}
	\sum_{n=1}^{d_k} \dbinom{N_k}{n}(1-S_k)^nS_k^{N_k-n}\leq\beta\,,
\end{equation*}
where $d_k$ is the number of optimization variables. The explicit lower bound of $N_k$ can be further derived as \cite{Calafiore09}:   
\begin{equation}\label{size}
	N_k\geq\frac{2}{1-S_k}\left(\ln \frac{1}{\beta}+d_k\right).
\end{equation}
To further reduce the conservatism of the scenario-based approach, a sample removal approach is proposed in \cite{Calafiore2011} and several variants are proposed. Their use here is beyond the scope of this paper.

For each scenario $i$, $n_\mathrm{c}$ linear constraints are generated in \eqref{deterministic}. It is clear that $n_\mathrm{c}$ is equal to the number of rows of matrix $H^{(i)}$. We additionally observe that for the constraint of index $j\in \mathbb I_{[1,n_\mathrm{c}]} $ in \eqref{deterministic}, the following inequality holds:
\begin{equation*}
	[H^{(i)}]_j+ [\mathcal C]_j\bar{\vect s}_{k+1,\ldots,N} \leq \max_{q\in\mathbb I_{[1,N_k]}}\,[H^{(q)}]_j + [\mathcal C]_j\bar{\vect s}_{k+1,\ldots,N}
\end{equation*}
for all $i\in\mathbb I_{[1,N_k]}$, where $[\star]_j$ denotes the $j^{\mathrm{th}}$ row of the matrix $\star$. Note that this inequality is tight, i.e., for all constraint of index $j$ there always exists at least one sample of index $i$ that ensures the above inequality tight. Hence $j^{\mathrm{th}}$ constraint is satisfied for all realizations $i$ if they are satisfied for the one having the largest $[H^{(i)}]_j$.

Let us label: $$\mathcal I_j = \max_{i\in\mathbb I_{[1,N_k]}}\,[H^{(i)}]_j,\quad\forall j\in \mathbb I_{[1,n_\mathrm{c}]}.$$  then we have that constraint \eqref{deterministic} is equivalent to the following constraints:
\[\mathcal I_j + [\mathcal C]_j\bar{\vect s}_{k+1,\ldots,N} \leq 0, \quad \forall j\in \mathbb I_{[1,n_\mathrm{c}]}.\]
Note that calculating $\mathcal I_j,$ for all $j\in \mathbb I_{[1,n_\mathrm{c}]}$, requires only
$n_\mathrm{c}$ (vector) maximum operations that are easy to implement and computationally efficient.


Now,  \eqref{eq:LinearSMPC:OCP} is equivalent to the following QP:
\begin{subequations}
	\label{eq:MWPS:QP}
	\begin{align}
		\min_{\bar{\vect u}_{k,\ldots,N-1}}&\quad \bar{\vect s}_N^{\mathrm{T}} Q_N \bar{\vect s}_N + \sum_{t=k}^{N-1}\left(\bar{\vect s}_t^{\mathrm{T}} Q \bar{\vect s}_t +  \bar{\vect u}_t^{\mathrm{T}} R \bar{\vect u}_t\right)\\
		\mathrm{s.t.}\,\,\,\,\, 
		&\quad \bar{ \vect s}_k=\vect s_k\\
		&\quad \bar{\vect s}_{t+1} = A\bar{\vect s}_t + B\bar{\vect u}_t + \bar{\vect w}_t,\,\, \forall t\in \mathbb I_{[k,N-1]} \\
		&\quad \mathcal I_j + [\mathcal C]_j\bar{\vect s}_{k+1,\ldots,N} \leq 0, \quad \forall j\in \mathbb I_{[1,n_\mathrm{c}]}
	\end{align}
\end{subequations}
yielding a regular QP of the same complexity as a normal linear MPC. 

A systematic overview of the proposed scenario-based mission-wide linear SMPC scheme is summarized in Algorithm \ref{alg:SMPC}.
\begin{algorithm}
	\caption{linear SMPC with MWPS constraints}
	\label{alg:SMPC}
	\DontPrintSemicolon
	\textbf{Initialization:} $S_0$, $\gamma_{1,\ldots,N-1}$, initial state $\vect s_0$ \\
	
	\For{$k=0:N-1$,}{
		1) \If{$k\geq 1$}{ 
			Evaluate $S_k$ in \eqref{Sk} through Monte Carlo simulation}
		2) Generate $N_k$ scenarios according to \eqref{size}\\
		3) Get the solution $\bar{\vect u}^*_{k,\ldots,N-1}$ by solving \eqref{eq:MWPS:QP}\\
		4) Send $\bar{\vect u}_k^*$ to the actual system and update state: $\vect s_{k+1} = A\vect s_k+B\bar{\vect u}_k^*+\vect w_k$
		
	}
\end{algorithm}

\subsection{Numerical Case Study}
We consider the linear system \eqref{eq:LinearSystem} with
\begin{align*}
	A = \matr{cc}{1 &1\\ 0 &1},\quad B = \matr{c}{0.5\\1}
\end{align*} 
and the uncertainty is assumed to have a Gauss distribution $$\vect w_k \sim \mathcal N (0,0.04\cdot I).$$
The safe (constraint) set $\mathbb S$ \eqref{eq:safeSet} is given by 
matrices
\begin{equation*}
	C=\matr{cc}{1 &0 \\0 &1\\-1 &0 \\0 &-1}\,,\quad \vect c = \matr{c}{-2 \\ -2 \\-10 \\-2}\,.
\end{equation*}
The matrices $Q = I$, $R = 0.1$, and  
\begin{equation*}
	K= [-0.6167, -1.2703]\,,\quad Q_N=\matr{cc}{2.0599 &0.5916 \\0.5916 &1.4228}
\end{equation*}
are computed from the corresponding LQR solution.

We select $N = 11$, $S_0=0.98$ and $\gamma_{1,\ldots,10}=0.99$, resulting in $S = \prod_{k=1}^{10}\gamma_k S_0 = 0.8863$. The number $N_k$ of disturbance sample is selected from \eqref{size}. The bound $S_k$ given by \eqref{Sk} is evaluated from Monte Carlo simulation and $\beta=10^{-6}$. In the simulations, we observed that $S_k\approx0.99$ for all $k\in\mathbb I_{[1,10]}$. This is due to $N_k$ calculated from \eqref{size} is conservative, such that the remaining MWPS at time $k$ achieved by the previous policy sequence $\{\vect \pi^{k-1}_k,\ldots,\vect \pi^{k-1}_{N-1}\}$, is much higher than that is actually required.

A Monte Carlo simulation that simulates $10^5$ missions shows that the resulting ratio of mission success is $99.88\%$. This result is larger than $S=88.63\%$. The reason for this discrepancy is that the scenario-based method adopted is conservative. Fig. 2 shows the state trajectories of $10^3$ missions. 

\begin{figure}[htb]
	\centering
	\begin{overpic}[width=0.49\textwidth]{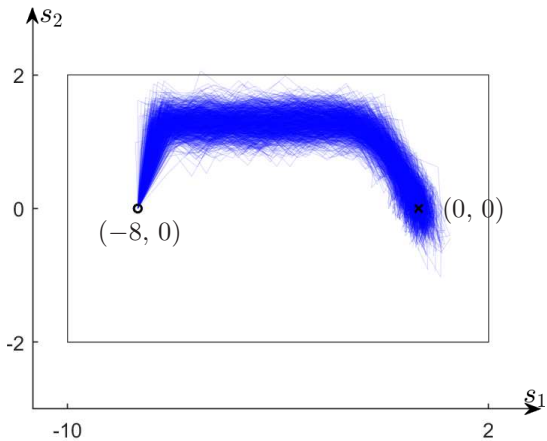}
		\put(75.5,37.5){$(0,\, 0)$}
		\put(23,34){$(-8,\, 0)$}

	\end{overpic}
	\caption{State trajectories plot obtained by running $10^{3}$ number of missions starting from the initial sate $\vect s_0=[-8,\,0]^\top$. The reference point is  $[0,\,0]^\top$. The rectangular area depicts the safe set $\mathbb S$.}
	\label{fig:trajectory}
\end{figure}

\section{Conclusions and Future Work}\label{sec:conclusion}
We investigated optimal policies satisfying Mission-Wide Probability of Safety constraints, i.e. constraints imposing the safety of a system over an entire mission. This is in contrast with classical stochastic MPC, where safety constraints are imposed independently at every time stage. We show that recursive feasibility holds in the expected value sense for the concept of Mission-Wide Probability of Safety, opening a simple and practically meaningful concept of recursive feasibility for stochastic MPC. Optimal control with mission-wide probabilistic constraints is challenging. However, a computationally efficient scenario-based approach is proposed to solve this issue for linear stochastic problems. For the sake of brevity, a shrinking-horizon approach was presented in this paper. The scenario-based approach proposed here relies on classical Monte-Carlo sampling. More advanced methods will be developed in the future for the proposed method.
 



\bibliography{references}
\bibliographystyle{ieeetr} 

\end{document}